\documentclass[12pt]{amsart}
\usepackage{graphics}
\usepackage{amsthm}
\usepackage{amssymb}
\usepackage{amsmath}
\usepackage{amsfonts} 
\input{psfig.sty}

\newcommand{\baseRing}[1]{\ensuremath{\mathbb{#1}}}
\newcommand{\Z}{\baseRing{Z}}
\newcommand{\C}{\baseRing{C}}
\newcommand{\N}{\baseRing{N}}
\newcommand{\R}{\baseRing{R}}

\theoremstyle{plain}
\newtheorem{theorem}{Theorem}[section]
\newtheorem{lemma}[theorem]{Lemma}

\newtheorem{proposition}[theorem]{Proposition}

\newtheorem{conjecture}[theorem]{Conjecture}

\theoremstyle{definition}
\newtheorem{definition}[theorem]{Definition}
\newtheorem{observation}[theorem]{Observation}
\newtheorem{remark}[theorem]{Remark}

\def\Av{A \cdot (\eta+\alpha)}
\def\rank{ {\rm rank}\, }
\def\vol{ {\rm vol}\, }

\def\conv{ {\rm conv} }
\def\nsupp{ {\rm nsupp} }
\def\ini{ {\rm in}\, }

\def\ker{ {\rm ker}\, }
\def\coker{ {\rm coker}\, }
\def\re{ {\rm Re}\, }
\def\Span{ {\rm Span}\, }

\def\im{ {\rm im}\, }

\title[Exceptional parameters for $A$-hypergeometric systems]%
{Exceptional parameters for \\ generic $A$-hypergeometric systems}

\author{Laura Felicia Matusevich}
\thanks{The author was partially supported by a Julia B. Robinson 
fellowship at UC Berkeley.}
\address{University of California at Berkeley}
\email{laura@math.berkeley.edu}
\date{}

\begin{document}

\begin{abstract}
The holonomic rank of an $A$-hypergeometric system
$H_A(\beta)$ is conjectured to be independent of the parameter
vector $\beta$ if and only if the toric ideal $I_A$ is Cohen
Macaulay. We prove this conjecture in the case that $I_A$ is generic 
by explicitly constructing more 
than $\vol(A)$ many linearly independent hypergeometric functions
for parameters $\beta$ coming from embedded primes
of certain initial ideals of $I_A$. 
\end{abstract}

\maketitle

\section{Introduction}
\label{intro}

$A$-hypergeometric systems
are systems of linear partial differential equations with polynomial 
coefficients
that can be built out of a toric ideal and a parameter vector. Homogeneous
toric ideals are themselves built out of combinatorial data: $n$ 
distinct integer points
lying in a hyperplane off the origin in $d$-dimensional space. 
We may assume that these points are
the columns of a $d\times n$ integer matrix whose first row is made up of
ones.

\begin{definition}
A $d \times n$ matrix $A$ 
whose columns are distinct elements of $\{ 1 \} \times \Z^{d-1}$ and 
generate $\Z^d$ as a lattice 
is said to be {\bf homogeneous}.
We set $m=n-d$.
\end{definition}

As we have already mentioned, we will think of $A$ as a point 
configuration in $d$-space. These points (the columns of $A$)
will be called $a_1,\dots ,a_n$.
The convex hull of $\{ a_1,\dots ,a_n\}$, $\conv(A)$, is a $(d-1)$-dimensional
polytope, whose normalized volume we denote by $\vol(A)$.

\begin{definition}
Given a homogeneous matrix $A$, the {\bf toric ideal $I_A$} is the ideal
of the polynomial ring $\C[\partial_1, \dots ,\partial_n]$ given by
\[ I_A = \langle \partial^u-\partial^v : u,v \in \N^n, A\cdot u=A\cdot v \rangle \; .\]
Here we use multi-index notation $\partial^u = \partial_1^{u_1} \cdots 
\partial_n^{u_n}$.
\end{definition} 

$A$-hypergeometric systems are non-commutative objects, they are left ideals in
the Weyl algebra $D$. 
The Weyl algebra is the quotient of the
free associative algebra with generators 
$x_1,\dots ,x_n,\partial_1,\dots ,\partial_n$, modulo the relations
\[ x_ix_j=x_jx_j\; ,\; \partial_i\partial_j=\partial_j\partial_i \; ,\; 
\partial_i x_j = x_j \partial_i + \delta_{ij} \; , \; 1 \leq i,j \leq n \; ,\]
where $\delta_{ij}$ is the Kronecker delta.

\begin{definition}
Given a homogeneous matrix $A$ and a vector $\beta \in \C^d$, 
the {\bf $A$-hypergeome\-tric system with parameter vector $\beta$} is the
left ideal in the Weyl algebra generated by $I_A$ and the homogeneity 
operators:
\[ (\sum_{j=1}^n a_{ij}x_j \partial_j) - \beta_i \;\; , \;\; 
1 \leq i \leq d\; . \]
\end{definition}

Hypergeometric systems were introduced in the late eighties by Gel$'$fand, 
Kapranov and
Zelevinsky (see, for instance, \cite{GKZ}). 
As we define them here, they are regular holonomic.
This means two things. First, that the {\em holonomic rank} of $H_A(\beta)$, 
that is, the dimension of the space of holomorphic solutions of 
$H_A(\beta)$ around a 
nonsingular point, is well defined, and second, that the
solutions of $H_A(\beta)$ can be represented as 
power series with logarithms. The holonomic rank of $H_A(\beta)$ is denoted
$\rank(H_A(\beta))$. One of the first results shown by Gel$'$fand, 
Kapranov and Zelevinsky
about $A$-hypergeometric systems is that, when the underlying toric ideal $I_A$
is Cohen Macaulay, $\rank(H_A(\beta))=\vol(A)$ for all $\beta \in \C^d$. 
A proof of this theorem can be found in \cite[Section 4.3]{SST}. 
This equality can fail if $I_A$ is not Cohen Macaulay.

However,
as is shown in \cite{adolphson}, \cite[Theorem 3.5.1,Equation 4.3]{SST}, 
if we drop the Cohen Macaulayness hypothesis, $\rank(H_A(\beta)) \geq \vol(A)$
for all $\beta \in \C^d$, and equality holds for generic $\beta$.

We define the {\em exceptional set} of $A$ to be:
\[ {\mathcal{E}}(A) := \{ \beta \in \C^d : \rank(H_A(\beta))>\vol(A) \} \; .\]
Elements of ${\mathcal{E}}(A)$ are called {\em exceptional parameters}.
The following conjecture, due to Bernd Sturmfels, relates the
existence of exceptional parameters and the Cohen Macaulayness of $I_A$.

\begin{conjecture}
\label{conjecture}
The exceptional set of $A$ is empty if and only if $I_A$ is Cohen Macaulay.
\end{conjecture}

The ``if" part of this conjecture is the aforementioned result by Gel$'$fand, 
Kapranov and Zelevinsky. As for the ``only if" part, it has been
proved in various special cases.
The $d=2$ case of Conjecture \ref{conjecture} is a result of
Cattani, D'Andrea and Dickenstein (see \cite{monomialcurves}).
The $n-d=2$ case was proved by the author in \cite{codim2}. 
Finally, Saito has shown that Conjecture \ref{conjecture} 
is true when the convex hull of the
configuration $A$ is a simplex (see \cite{mutsumi2}).

In this article, we study Conjecture \ref{conjecture} in the case that
$I_A$ is a generic toric ideal.
\begin{definition}
A lattice ideal is {\bf generic} if it has a minimal generating
set of binomials with full support.
\end{definition}
Generic lattice ideals were introduced by Irena Peeva and Bernd Sturmfels
in \cite{genericlatticeideals}. One of the results in that article
is that ``most'' toric ideals are generic. 

The following is the main result in this article.

\begin{theorem}
\label{thm:maintheorem}
Let $I_A$ be a generic non Cohen Macaulay toric ideal. Then
the exceptional set of $A$, ${\mathcal{E}}(A)$ contains an
affine space of dimension $d-2$. In particular, 
${\mathcal{E}}(A) \neq \emptyset$.
\end{theorem}

In order to prove this theorem, we will need a characterization
of Cohen Macaulayness for generic toric ideals. 
We will use Lemma \ref{lemma:embeddedprimes}, that relates this
notion to the existence of embedded primes of certain initial ideals 
of our toric ideal.

This article is organized as follows. Sections \ref{sec:toricbackground}
and \ref{stdpairsandlogfreeseries} contain background material
on toric algebra and hypergeometric functions respectively.
In Sections \ref{sec:proof1}, \ref{sec:proofkernel},
\ref{sec:proofcokernel}, and \ref{sec:proofjump}, we prove Theorem
\ref{thm:maintheorem}. Section \ref{sec:example}
contains a completely worked out example, and in Section 
\ref{sec:comprehensive-gb} we show that, under no hypotheses on $A$,
the exceptional set ${\mathcal{E}}(A)$ is Zariski constructible.


\section{Standard pairs of initial ideals of toric ideals}
\label{sec:toricbackground}

In this section we introduce the notion of standard pairs. These objects
play a fundamental role in the study of the associated primes of monomial 
ideals: if $(\partial^{\eta},\sigma)$ is a standard pair of
a monomial ideal $M \subset \C[\partial_1,\dots,\partial_n]$, then $\langle
\partial_j : j \not \in \sigma \rangle$ is an associated prime
of $M$. Moreover, all associated primes of $M$ arise this way.

In the special case of initial ideals of toric ideals, standard pairs
admit a polyhedral description (Theorem \ref{thm:polytopeforstdpairs}).

\begin{definition}
Let $M$ be a monomial ideal of $\C[\partial_1,\dots,\partial_n]$. A 
{\bf standard 
pair} of $M$ is a pair $(\partial^{\eta}, \sigma)$, where $\eta \in \N^n$ and
$\sigma \subset \{ 1, \dots, n \}$ subject to the following three conditions:
\begin{enumerate}
\item $\eta_i = 0$ for $i \in \sigma$;
\item For all choices of integers $\mu_i \geq 0$, $i \in \sigma$,
the monomial
$\partial^{\eta} \cdot \prod_{i \in \sigma} \partial_i^{\mu_i} $ is not in $M$.
\item For all $l \not \in \sigma$, there exist $\mu_i \geq 0$,
$i \in \sigma \cup \{ l \}$, such that
$\partial^{\eta} \cdot \partial_l^{\mu_l} \cdot 
\prod_{i \in \sigma} \partial_i^{\mu_i} $
is in $M$.
\end{enumerate}
The set of standard pairs of $M$ is denoted $S(M)$. A standard pair 
$(\partial^{\eta}, \sigma)$ such that the ideal $\langle \partial_i
: i \not \in \sigma \rangle$ is a minimal associated prime of $M$
is called {\bf top-dimensional}.  Standard pairs
that are not top-dimensional are called {\bf embedded}.
\end{definition}

If $A$ is a homogeneous $d\times n$ matrix and $w\in \R^n$ is a generic weight
vector for $I_A$, that is, $\ini_{w}(I_A)$ is a monomial ideal, we can
study $S(\ini_{w}(I_A))$ using combinatorial techniques. For instance, a
standard pair of $\ini_w(I_A)$ is top-dimensional if and only if
the cardinality of $\{1, \dots ,n\} \backslash \sigma$ 
is equal to $m$ (see \cite[Corollary 2.9]{associatedprimes}).
The following polyhedral characterization of standard pairs
of initial ideals of toric ideals
is due to Serkan Ho\c{s}ten and Rekha Thomas 
(see \cite[Theorems 2.3, 2.5]{associatedprimes}). 
Choose a $\Z$-basis of $\ker_{\Z}(A)$ and form an $n\times m$ 
matrix $B=(b_{ij})$
whose columns are the vectors in this basis. This matrix $B$
is called a {\em Gale dual} of $A$.

\begin{theorem}
\label{thm:polytopeforstdpairs}
A pair $(\partial^{\eta}, \sigma)$, where $\eta \in \N^n$ and
$\eta_i = 0$ for $i \in \sigma$, is a standard pair of 
the monomial ideal $\ini_w(I_A)$
if and only if $0$ is the only lattice point in the polytope
\[ P_{\eta}^{\bar{\sigma}}:=
\{ y \in \R^m : (B\cdot y)_j \leq \eta_j \;, j \not \in \sigma\; ;
-(w)^t(B\cdot y) \leq 0 \}\; , \]
and all the inequalities $(B\cdot y)_j \leq \eta_j \; , j \not \in \sigma$,
are essential, that is, removing an inequality introduces a new lattice
point $z$ into the resulting polyhedron. We may assume that
$z$ is such that $-w^t(B\cdot z)$ is strictly negative.
\end{theorem}

The fact that $P^{\bar{\sigma}}_{\eta}$ is a polytope (and hence
a bounded set) has the following linear algebra consequence.

\begin{lemma}
\label{lemma:lind}
Let $(\partial^{\eta},\sigma)$ be a standard pair of 
$\ini_{w}(I_A)$, where $w$ is such that this is a monomial ideal.
Then the set 
$\{(b_{i1},\dots ,b_{im}) : i \not \in \sigma\}$
contains a linearly independent subset
of cardinality $m$.
\end{lemma}

\begin{proof}
By contradiction, suppose that no subset of cardinality $m$
of $\{ (b_{i1},\dots ,b_{im}): i \not \in \sigma\}$ is linearly independent. 
This means that the matrix whose rows are the rows of $B$ 
indexed by $i \not \in \sigma$ has rank strictly less than $m$.
Consequently, we
can find rational numbers $s_1,\dots,s_m$ not all zero, such that
$(B\cdot (s_1,\dots ,s_m)^t)_i=0$, for all $i \not \in \sigma$. 
But then at least half of the line $\{ \lambda (s_1,\dots,s_m)^t : 
\lambda \in \R\}$
is contained in $P_{\eta}^{\bar{\sigma}}$, contradicting that this
set is bounded.
\end{proof}

Associated primes of initial ideals of toric ideals
come in saturated chains. This is the content of the following
theorem, due to Ho\c{s}ten and Thomas \cite[Theorem 3.1]{associatedprimes}.

\begin{theorem}
\label{thm:chaintheorem}
Let $I_A$ be a toric ideal 
and $w\in \R^n$ a generic weight vector for $I_A$. If ${\mathfrak{p}}$ is an
embedded prime of $\ini_w(I_A)$, then ${\mathfrak{p}}$ contains an
associated prime ${\mathfrak{q}}$ of $\ini_w(I_A)$ such that 
$dim({\mathfrak{q}})=dim({\mathfrak{p}})+1$.
\end{theorem}

Finally, we include here a combinatorial lemma that will be useful later on.
I am grateful to Bernd Sturmfels, who provided this beautiful proof.

\begin{lemma}
\label{lemma:boundedregions}
Let $c^{1}, \dots , c^{m+2} \in \R^m$ and $k_1, \dots k_{m+2} \in \R$
and consider $P =  \{ z \in \R^m : c^{j}\cdot z  \leq k_j, j=1,\dots m+2 \}$.
Suppose that this set is nonempty.
There is a set $T \subset \{ 1, \dots , m+1 \}$ of cardinality $m$
such that, for each
$j \in T$, the set obtained from $P$ by reversing the inequality
$c^j \cdot z \leq k_j$ is unbounded. Moreover, the set
$\{ c^{j}: j \in T \}$ is linearly independent.
\end{lemma}

\begin{proof}


We use the technique of Gale duality, as introduced in 
\cite[Chapter 6]{leconpolytopes}. Consider the vector
configuration ${\mathcal{C}} = \{c_1,\dots ,c_{m+2} \} \subset \R^m$. Let
${\mathcal{B}} = \{b_1,\dots ,b_{m+2} \} \subset \R^2$ be its Gale dual 
configuration. 
We define ${\mathcal{N}}$ as the set of indices $j$ such that
reversing the $j$-th inequality in $P$ produces a bounded set.
Now, $i \in {\mathcal{N}}$ if an only if the
configuration ${\mathcal{C}} \backslash \{ c_i \}$ is totally
cyclic, that is, if and only if there exists a positive linear dependence
among the elements of ${\mathcal{C}} \backslash \{ c_i \}$.

Using results from \cite[Sections 6.3(b), 6.4]{leconpolytopes},
namely that the dual operation to deletion is contraction, and that
dual of a totally cyclic vector configuration is an acyclic vector
configuration, we see that $i \in {\mathcal{N}}$
if and only if ${\mathcal{B}}/b_i \subset \R$ 
is acyclic, that is, if and only if
there is a linear functional $c \in \R$ such that 
$c \cdot b > 0 $ for all $b \in {\mathcal{B}}/b_i$. But this can only
happen if $b_i$ is an extreme ray of ${\mathcal{B}}$, which means that
there exists $c' \in \R^2$ such that $c' \cdot b_i = 0$
and $c' \cdot b_j > 0$ for all $1 \leq j \leq m+2$, $j \neq i$.
Since a configuration in $\R^2$ has at most two extreme rays, we
conclude that the cardinality of ${\mathcal{N}}$
is at most $2$.
The last assertion of this lemma follows from
the same arguments used to prove Lemma \ref{lemma:lind}.
\end{proof}


\section{Canonical Hypergeometric Series}
\label{stdpairsandlogfreeseries}

In this section we review material about canonical logarithm-free
hypergeometric series. Our source is \cite[Sections 2.5, 3.1, 3.4, 4.1]{SST}.
We start by relating standard pairs to hypergeometric functions through
the concept of fake exponents.

\begin{definition}
\label{fakeexponents}
Let $A$ be a homogeneous $d \times n$ matrix, $\beta \in \C^d$, 
and $w \in \R^n$ a generic weight
vector for $I_A$. The set of 
{\bf fake exponents of $H_A(\beta)$ with respect to $w$}
is the vanishing set of the following zero dimensional ideal
of the (commutative) polynomial ring $\C[\theta_1,\dots ,\theta_n]$:
\[  \bigcap_{(\partial^{\eta},\sigma) \in S(\ini_w(I_A))} 
(\langle \theta_i-\eta_i : i \not \in \sigma \rangle + \langle 
A \cdot \theta - \beta 
\rangle ) \, , \]
where $\theta$ is the vector $(\theta_1,\dots ,\theta_n)^t$.
\end{definition}

If $A$ is homogeneous and $\beta \in \C^d$, the hypergeometric system
$H_A(\beta)$ is regular holonomic. This means that the solutions
of $H_A(\beta)$ can be written as power series with logarithms. 
Moreover, if $w$ is a generic weight vector for $I_A$, we can
use the techniques in \cite[Sections 2.5, 3.4, 4.1]{SST} to 
build a basis of the solution space of $H_A(\beta)$,
whose elements are called {\em canonical series} (with respect to $w$). 
This basis is defined by Proposition \ref{canonicalseries1}.

Notice that
we can extend the notion of term order to the ring of power series 
with logarithms as follows:
\begin{equation}
x^{\alpha} \log(x)^{\gamma} \leq x^{\alpha'} \log(x)^{\gamma'}
 \; \Longleftrightarrow 
\; \re(w\cdot \alpha) \leq \re(w\cdot \alpha') . \notag
\end{equation}
Here we mean $x^{\alpha}=x^{\alpha_1}\cdots x^{\alpha_n}$ for 
$\alpha \in \C^n$,
where $x_i^{\alpha_i}=\exp(\alpha_i \log(x_i))$; and 
$\log(x)^{\gamma}=\log(x_1)^{\gamma_1} 
\cdots \log(x_n)^{\gamma_n}$, for $\gamma \in \N^n$.
If we refine this ordering lexicographically, we can
define the initial term
of a power series with logarithms $\varphi$ (if it exists)
as $\ini_{\prec_w}(\varphi) = \;$the {\em least} term of $\varphi$
with respect to this ordering. By \cite[Proposition 2.5.2]{SST},
$\ini_{\prec_w}(\varphi)$ exists when $\varphi$ is a solution of $H_A(\beta)$
that converges in a certain region of $\C^n$ that depends on $w$.

\begin{proposition}
\label{canonicalseries1}
There exist a basis of the solution space of $H_A(\beta)$
such that
\begin{itemize} 
\item if
$\varphi$ is an element of this basis, 
there exists
a fake exponent $v$ of $H_A(\beta)$ with respect to $w$ such that
$\ini_{\prec_w}(\varphi)=x^v\log(x)^{\gamma}$ for some $\gamma \in \N^n$,
\item if $\varphi$ and $\psi$ are elements of this basis,
then $\ini_{\prec_w}(\varphi)$ is not a term appearing in $\psi$.
\end{itemize}
This basis is called the {\bf basis of canonical series with respect to $w$}.
\end{proposition}

\begin{proof}
This follows from Corollary 2.5.11, Corollary 3.1.6 and 
Lemma 4.1.3 in \cite{SST}.
\end{proof}

\begin{remark}
There might be fake exponents of $H_A(\beta)$ with respect to $w$
that do not give initial monomials of canonical series (see
Example 3.1.8 in \cite{SST}). Also the same fake exponent 
might give rise to two different canonical series. If this is the case,
then at least one of those canonical series will have logarithms.
\end{remark}

Using the results from \cite[Section 3.4]{SST} we can pinpoint exactly
which fake exponents give initial monomials of {\em logarithm-free}
canonical series solutions of $H_A(\beta)$.
To do this we need to introduce the following concepts.
The {\em negative support} of a vector $v \in \C^n$ is the
set:
\[ \nsupp(v):=\{i \in \{1,\dots ,n\} : v_i \in \Z_{<0} \}. \]
A vector $v \in \C^n$ has {\em minimal negative support} 
with respect to $A$ if
\begin{align}
& u \in \ker_{\Z}(A) \; \mbox{and} \; \nsupp(v-u) \subseteq 
  \nsupp(v) \; \mbox{imply} &  \notag \\
 &\nsupp(v-u)=\nsupp(v) \, . & \notag
\end{align}
In this case, let 
\[ N_v = \{ u \in \ker_{\Z}(A) : \nsupp(v-u)=\nsupp(v) \} \, , \]
and define the following formal power series:
\begin{equation}
\label{canonicalseries}
\phi_v = \sum_{-u \in N_v} \frac{[v]_{u_-}}{[u+v]_{u_+}} x^{v+u} \;\; ,
\end{equation}
where
\[ [v]_{u_-} = \prod_{i:u_i<0} \prod_{j=1}^{-u_i} (v_i-j+1) 
\;\;\;\; \mbox{and} 
\;\;\;\;
 [u+v]_{u_+} = \prod_{i:u_i>0} \prod_{j=1}^{u_i} (v_i+j) \, . \]

\begin{theorem}{\cite[Theorem 3.4.14, Corollary 3.4.15]{SST}}
Let $v \in \C^n$ be a fake exponent of $H_A(\beta)$ with minimal negative
support. Then the series $\phi_v$ defined in 
(\ref{canonicalseries}) is a canonical solution of the $A$-hypergeometric
system $H_A(\beta)$. In particular, $\phi_v$ converges in a region of $\C^n$. 
The set:
\[ \{ \phi_v : v \; \mbox{is a fake exponent with minimal negative support} \, \} \]
is a basis of the space of logarithm-free solutions of $H_A(\beta)$.
\end{theorem}

We now mention a way to distinguish which vectors with minimal negative
support are fake exponents.

\begin{proposition}
Let $v \in \C^n$ with minimal negative support such that $A\cdot v = \beta$.
Then $v$ is a (fake) exponent of $H_A(\beta)$ with respect to a weight 
vector $w$ if and only if 
\[ w \cdot v = \mbox{min} \, \{ w \cdot u : \nsupp(u)=\nsupp(v) \,
\mbox{and} \; u-v \in \ker_{\Z}(A) \} \]
and this minimum is attained uniquely.
\end{proposition}

At the moment, a detailed characterization such as we have for 
logarithm-free $A$-hypergeo\-metric series does not exist
for logarithmic $A$-hypergeometric functions. However, information
about logarithmic series will be necessary to prove Theorem
\ref{thm:maintheorem}. We reproduce two results about logarithmic
hypergeometric functions. The first is an observation from
\cite{codim2}.

\begin{observation}
\label{keyobservation}
Let $\psi$ be a solution of $H_A(A\cdot v)$.
This function is of the form:
\[ \psi = \sum 
c_{\alpha,\gamma}x^{\alpha} \log(x)^{\gamma} \, ,\]
where the sum runs over $\alpha$ such that $A\alpha=A\cdot v$, 
$\gamma \in \{0,1,\dots ,h-1\}^n$ for $h = \rank(H_A(A\cdot v))$.

The set ${\mathcal{S}}:=\{ \gamma \in [0,h-1]^n\cap\N^n : \exists 
\alpha \in \C^n \,
\mbox{such that} \; c_{\alpha,\gamma} \neq 0 \}$ 
is partially ordered with respect to:
\[ (\gamma_1,\ldots,\gamma_n) \leq (\gamma_1',\ldots,\gamma_n') 
\Longleftrightarrow \gamma_i \leq \gamma_i', \;
i=1,\ldots , n \, .\]

Denote by ${\mathcal{S}}_{\max}$ the set of maximal elements of 
${\mathcal{S}}$.
Let $\delta \in {\mathcal{S}}_{\max}$ and  $f_{\delta} = \sum_{\alpha \in \C^n}
c_{\alpha,\delta}x^{\alpha}$. Write
\[ \psi = \psi_{\delta} + \log(x)^{\delta} \, f_{\delta} \, ,\]
so that the logarithmic terms in $\psi_{\delta}$ are either less than
or incomparable to $\delta$.
If $P$ is a differential operator that annihilates $\psi$, we have:
\[0 = P \psi = P\psi_{\delta}+\log(x)^{\delta} Pf_{\delta}+
\mbox{terms whose $\log$ factor is lower than
$\delta$} .\]
Since $P\psi_{\delta}$ is a sum of terms whose $\log$ factor is 
either lower than $\delta$
or incomparable to $\delta$, we conclude that $Pf_{\delta}$ must be zero. 
This implies that
$f_{\delta}$ 
is a logarithm-free $A$-hypergeometric function of degree $A\cdot v$.
{\em Moreover, if $\partial_1 \psi$ is logarithm-free, then
$\partial_1 f_{\delta}$ must vanish.}
\end{observation}

The following result is due to Saito (see \cite{mutsumi2}).

\begin{proposition}
\label{propo:formoflogs}
Let $\varphi = \sum x^u g_u(\log(x))$ be a solution of $H_A(\beta)$. Then
the polynomials $g_u$ are of the form:
\[ g_u (t_1,\dots,t_n) = c_0 + c^{(1,1)}\cdot (t_1,\dots,t_n)+\cdots
+ \prod_{i=1}^l c^{(l,i)} \cdot (t_1,\dots,t_n) , \]
where $c_0 \in \C$, the vectors $c^{(i,j)}$ belong to the kernel of $A$,
and 
\[c^{(i,j)} \cdot (t_1,\dots ,t_n) = c^{(i_j)}_1t_1 + \cdots+
c^{(i,j)}_n t_n \; .\; \]
Moreover, if $g_u \neq 0$, then $A\cdot u = \beta$.
\end{proposition}

Finally, the techniques from \cite[Section 3.5]{SST} imply the following 
proposition.

\begin{proposition}
\label{propo:closure}
Let $\beta, \beta' \in \C^n$ and suppose that
$\rank(H_A(\beta+\epsilon \beta')) \geq t$ for $0< \epsilon < \epsilon_0$.
Then $\rank(H_A(\beta)) \geq t$.
\end{proposition}


\section{Constructing exceptional parameters}
\label{sec:proof1}

In this section we start working towards
the proof of Theorem \ref{thm:maintheorem}.
To do this, we will use the following characterization
of the Cohen Macaulay property for generic toric ideals
via associated primes of reverse 
lexicographic initial ideals of $I_A$.

\begin{lemma}
\label{lemma:embeddedprimes}
Let $I_A$ be a generic toric ideal. Then $\ini_{-e_i}(I_A)$
is a monomial ideal for all $1 \leq i \leq n$.
Moreover, for a generic toric ideal, the following are equivalent:
\begin{enumerate}
\item \label{it1} $I_A$ is Cohen Macaulay,
\item \label{it2} For all $i$, $\ini_{-e_i}(I_A)$ is free of embedded primes,
\item \label{it3} For all $i$, $\ini_{-e_i}(I_A)$ is Cohen Macaulay.
\end{enumerate}
\end{lemma}

\begin{proof}

The first assertion follows from \cite[Lemma 8.4]{syzygiesofcodim2}.
To see that (\ref{it1}) implies (\ref{it2}), notice that,
since $I_A$ is a Cohen Macaulay prime ideal , 
$I_A+\langle \partial_i \rangle$ is Cohen Macaulay, and thus free of embedded
primes, for all $1 \leq i \leq n$. As a consequence, for each $i$,
the ideal $\ini_{-e_i}(I_A)$, which is the ideal of
$\C[\partial_1,\dots \partial_n]$ generated by 
$(I_A+\langle \partial_i \rangle) \cap \C[\partial_1,\dots,\hat{\partial_i},
\dots,\partial_n]$, is free of embedded primes.
The second implication
follows from Theorems 2.5 and 3.1 in \cite{genericandcogeneric}. 
Finally, $\ini_{-e_i}(I_A)$ being
Cohen Macaulay implies that $I_A$ is Cohen Macaulay.
\end{proof}

\begin{remark}
From now on, we assume that $I_A$ is a generic non Cohen Macaulay
toric ideal.
\end{remark}

Consider $\ini_{-e_1}(I_A)$. By Lemma \ref{lemma:embeddedprimes}, 
this initial ideal
has an embedded prime. But since $\ini_{-e_1}(I_A)$ is a monomial ideal,
Theorem \ref{thm:chaintheorem} implies that  $\ini_{-e_1}(I_A)$
has an embedded prime of the form 
$\langle \partial_j : j \not \in \{ 1 \} \cup \tau \rangle$, where 
$\tau \subset \{ 2, \dots , n \}$ has cardinality $d-2$. By interchanging
the columns of $A$ we may assume that $\tau = \{ m+3,\dots,n\}$.
Since standard pairs carry all the information about the associated
primes of monomial ideals, we conclude that $\ini_{-e_1}(I_A)$
has an embedded standard pair $(\partial^{\nu},\{1\} \cup \tau)$.

\begin{observation}
\label{obs:changeofweight}
We may assume that the set obtained from 
$P^{\overline{\{1 \} \cup \tau}}_{\nu}$ by reversing the
weight inequality $(B\cdot z)_1 \leq 0$ is bounded, or that
the hyperplane $\{ z: (B\cdot z)_r= 0\}$
coincides with the hyperplane $\{ z: (B\cdot z)_1= 0\}$
for some $r \not \in \{ 1 \} \cup \tau$. In the latter case, the
first row of $B$ is a negative multiple of the $r$-th row of $B$, and
$\nu_r>0$.
\end{observation}

\begin{proof}

Suppose that the hyperplanes $\{ z: (B\cdot z)_i =0\}$,
for $i \not \in \{ 1 \} \cup \tau$, are pairwise distinct.
To check the  assertion, first notice that, since 
$\langle \partial_j : j \not \in \{ 1 \} \cup \tau \rangle$ is an associated
prime of $\ini_{-e_1}(I_A)$, it is contained in a minimal prime
$\langle \partial_j : j \not \in \{ 1, l \} \cup \tau \rangle$ of
this initial ideal. Here $l \in \{2,\dots,m+2\}$.
Thus, the set obtained from $P^{\overline{\{1 \} \cup \tau}}_{\nu}$
by removing the inequality $(B\cdot z)_l \leq \nu_l$ must be a simplex.
Now, by Theorem \ref{thm:polytopeforstdpairs}, any
lattice point $z\neq 0$ in this simplex satisfies $(B\cdot z)_l > \nu_l$,
and such lattice points exist. Pick $z \in \Z^m \backslash \{ 0 \}$ in
our simplex such that $(B\cdot z)_l$ is minimal. Notice that 
the lattice point $z$ is unique, since the weight vector $-e_l$ is generic.
This follows from results in \cite[Section 5]{intprog}.

Now let $\partial^{\mu} = 
\partial_1^{-(B\cdot z)_1-1}\prod_{i \in \{2,\dots,\hat{l},
\dots,m+2\}} \partial_i^{\nu_i-(B\cdot z)_i}$ and consider the pair
$(\partial^{\mu}, \{ l \} \cup \tau)$. It is easy to check that
this pair satisfies the conditions of Theorem
\ref{thm:polytopeforstdpairs} for the weight vector $-e_l$.
Thus, interchanging the first and $l$-th columns of $A$, and
replacing $\nu$ by $\mu$, we obtain a standard pair as we desired.

Finally, if the hyperplanes $\{ z : (B\cdot z)_r = 0 \}$
and $\{ z : (B\cdot z)_s = 0 \}$ coincide, for $r,s \not \in \{ 1 \} \cup 
\tau$, we use the above argument to change the weight vector $-e_1$
by the weight vector $-e_s$, and assume that $\{ z : (B\cdot z)_r = 0 \}$
and $\{ z : (B\cdot z)_1 = 0 \}$ are parallel.
Of course, the first row of $B$ is a multiple of the $r$-th row of $B$.
If it were a positive multiple, then for each $v \in \ker_{\Z}(A)$,
$\ini_{-e_1}(\partial^{v_+}-\partial^{v_-})$ does not contain the
variable $\partial_r$, which implies $r \in \tau$, a contradiction.
Thus the first and $r$-th rows of $B$ are negative multiples of each other.
To see that $\nu_r > 0$, notice that $\nu_r=0$ would contradict
the last assertion of Theorem \ref{thm:polytopeforstdpairs}.
\end{proof}

We are now ready to start the construction of our candidates for exceptional
parameters. 
Pick an embedded standard pair 
$(\partial^{\nu}, \{ 1 \} \cup \tau)$,
and consider the set 
\[ \left\{ (\partial^{\mu},\{ 1 \} \cup \tau) : \begin{array}{c}
(\partial^{\mu},\{ 1 \} \cup \tau) \in S(\ini_{-e_1}(I_A)), \\
\exists y_{\mu} \in \C^m
\; \mbox{such that} \; \mu_i = \nu_i-(B\cdot y_{\mu})_i \; \\ \mbox{for}\;  
2 \leq i \leq m+2, \; \mbox{and}\;0>(B\cdot y_{\mu})_1 \in \Z \end{array}
\right\}.\]
If this set is nonempty, select a standard pair $(\partial^{\eta},
\{ 1 \} \cup \tau)$ such that $(B\cdot y_{\eta})_1$ is maximal.
Otherwise, rename $\nu$ to $\eta$. This choice implies the following
fact.

\begin{observation}
\label{obs:noembeddedpairs}
If $(\partial^{\mu},\{ 1 \} \cup \tau)$ is a standard pair
of $\ini_{-e_1}(I_A)$ such that there exists $y \in \C^m$
that satisfies
\begin{enumerate}
\item $\mu_i = \eta_i - (B\cdot y)_i$ for $2 \leq i \leq m+2$, and
\item $(B\cdot y)_1 \in \Z$,
\end{enumerate}
then $(B\cdot y)_1 \leq  0$.
\end{observation}

Choose generic numbers $\alpha_i$ for $i \in \tau$. Let
$\alpha = \sum_{i \in \tau} \alpha_i e_i$. We look at the fake
exponents of $H_A(\Av)$ with respect to the weight vector $-e_1$.
Since the $\alpha_i$ are generic, we may assume that if $u$ is a 
fake exponent and $u_1=0$, then $u_i \in \Z$ implies 
$i \not \in \tau$. In particular, the numbers $\alpha_i$
are non integers. We define two sets:
\[ F := \bigg\{ \begin{array}{c}
\mbox{\rm fake exponents of $H_A(\Av)$ with respect to $-e_1$} \\
\mbox{\rm that have minimal negative support} \end{array} \bigg\} \]
and 
\[ K := \{ u \in F : u_1 = 0 \}. \]

\begin{proposition}
\label{propo:case1holds}
The following condition holds:
\begin{equation}
\label{cond:case1}
\mbox{If $u \in K$, $v \in F$, and $v-u \in \Z^m$, then $1 \not \in
\nsupp(v)$.}
\end{equation}
\end{proposition}

\begin{proof}

Suppose that Condition \ref{cond:case1} does not hold.
Then there are $u \in K$, $v \in F$ such that
$v-u\in \Z^m$, and $1 \in \nsupp(v)$.
Since both $u$ and $v$ have minimal negative support (and their
only integer coordinates are indexed by $i \not \in \tau$ by
the choice of $\alpha$) we conclude that
$\nsupp(u)=\{ s \}$ for some $ 2 \leq s \leq m+2$, and
$\nsupp(v)=\{ 1 \}$. 

Suppose that the standard pair that gives rise to $v$ is top-dimensional.
Then it is of the form $(\partial^{\mu},\{1,l\}\cup \tau)$
for some $2 \leq l \leq m+2$.
Since $u_s<0$, the standard pair corresponding to $u$ is of the form
$(\partial^{\nu},\{ 1,s \} \cup \tau)$.

Suppose $l=s$.
Remember that we have $z \in \Z^m$ such that $v-B\cdot z = u$. 
Then $(B \cdot z)_j \leq v_j=\mu_j$ for $j \not \in \{1,l\} \cup \tau$,
and $(B\cdot z)_1 < 0$. Then $z \in P^{\overline{\{1,l\}\cup \tau}}_{\mu}\cap
\Z^m = \{ 0 \}$, a contradiction. Thus $l \neq s$.

We consider two cases, as in the proof of Observation 
\ref{obs:changeofweight}. If two of the hyperplanes
$\{ y : (B\cdot y)_i=0 \}$, $i \not \in \tau$, coincide, then
$\{ y : (B\cdot y)_r=0 \}$ is equal to $\{ y : (B\cdot y)_1=0 \}$
for some $2 \leq r \leq m+2$. This implies that
$l=s=r$, because otherwise the sets $P^{\overline{\{1,s\}\cup\tau}}_{\nu}$
and $P^{\overline{\{1,l\}\cup\tau}}_{\mu}$ are unbounded. This contradicts
the previous paragraph.

Now suppose that the hyperplanes $\{ y : (B\cdot y)_i=0 \}$, $i \not \in \tau$,
are pairwise distinct. Since $P^{\overline{\{1,s\}\cup\tau}}_{\nu}$ and
$P^{\overline{\{1,l\}\cup\tau}}_{\mu}$ are simplices,
reversing the $s$-th or the $l$-th inequality in our original
polytope $P^{\overline{\{1 \} \cup \tau}}_{\eta}$ yields bounded sets.
By Observation \ref{obs:changeofweight}, if we reverse the inequality
$(B\cdot z)_1 \leq 0$ we also get a bounded set. In view of 
Lemma \ref{lemma:boundedregions}, we derive a contradiction.

In conclusion, the standard pair corresponding to $v$ cannot be 
top-dimensional, so it must be of the form 
$(\partial^{\mu},\{1 \} \cup \tau)$, by our choice of $\alpha$.
However, by Observation \ref{obs:noembeddedpairs}, a fake exponent 
of $H_A(\Av)$ coming from such a standard pair cannot have a
negative integer first coordinate.
This contradiction concludes the proof.
\end{proof}

In what follows, we will prove the following version
of Theorem \ref{thm:maintheorem}.

\begin{theorem}
\label{thm:restricted}
Suppose $\ini_{-e_1}(I_A)$ is a monomial ideal with an embedded
standard pair $(\partial^{\eta},\{1\} \cup \tau)$ as in Observations \ref{obs:changeofweight} and \ref{obs:noembeddedpairs}, where
$\tau=\{m+3,\dots ,n\}$.
Choose generic numbers $\alpha_i$ as above. Then 
\[ \beta := A\cdot (\eta+\alpha -e_1) \in {\mathcal{E}}(A). \]
\end{theorem}

\begin{remark}
{\em (Theorem \ref{thm:maintheorem} follows from Theorem 
\ref{thm:restricted}.)}
Notice that, once Theorem \ref{thm:restricted} is proved, 
Proposition \ref{propo:closure} 
will lift the assumption that
the numbers $\alpha_i$ are generic. Thus we will have produced
a $(d-2)$-dimensional affine space contained in ${\mathcal{E}}(A)$,
and the proof of Theorem \ref{thm:maintheorem} will be complete.
\end{remark}

\begin{remark}
From now on, we work under the hypotheses and
notation of Theorem \ref{thm:restricted}.
\end{remark}

We want to show that $\rank(H_A(\beta)) > \vol(A)$. One way to do this
is to show that $\rank(H_A(\beta)) > \rank(H_A(\Av))$, since
$\rank(H_A(\Av)) \geq \vol(A)$. The tool to compare these two
ranks is the $D$-module map (see \cite[Section 4.5]{SST})
\[ D/H_A(\beta) \longrightarrow D/H_A(\Av) \]
given by right multiplication by the operator $\partial_1$.
This induces a vector space homomorphism
between the solution spaces of $H_A(\Av)$ and $H_A(\beta)$:
if $\varphi$ is a solution of $H_A(\Av)$, then
$\partial_1 \varphi$ (the derivative of $\varphi$ with respect to
the variable $x_1$) is a solution of $H_A(\beta)$. It is this
vector space map that we want to study. More precisely, we want
information about the dimension of its kernel and cokernel.

\section{The kernel of the map $\partial_1$}
\label{sec:proofkernel}

We start our analysis of the map $\partial_1$ between
the solution spaces of $H_A(\Av)$ and $H_A(\beta)$ by describing
its kernel. The following proposition is the
first step in this direction.

\begin{lemma}
\label{lemma:monomial}
If $u \in K$, then the canonical series corresponding to $u$, $\phi_u$
equals the monomial $x^u$. Consequently, $\partial_1 \phi_u=0$.
Conversely, if $\phi_u$ is a logarithm-free canonical series such that
$\partial_1 \phi_u = 0$, then $u \in K$.
\end{lemma}

\begin{proof}

To see that $\phi_u$ is a monomial, it is enough to show that $N_u = \{ 0 \}$.
Remember that $N_u = \{ B\cdot z : z \in \Z^m \; \mbox{and} \; 
\nsupp(u-B\cdot z) = \nsupp(u) \}$.
Since $u$ is a fake exponent with respect to $-e_1$, there is a 
standard pair $(\partial^{\mu}, \{ 1 \} \cup \sigma)$ of
$\ini_{-e_1}(I_A)$ such that $u_i = \mu_i \in \N$ for 
$ i \not \in \{1 \} \cup \sigma$. Pick $z \in \Z^m$ such that
$B\cdot z \in N_u$. Then, since $\nsupp(u-B\cdot z) = \nsupp(u)$,
$(B\cdot z)_i \leq u_i=\mu_i$ for $i \not \in \{1\} \cup \sigma$, and
$(B\cdot z)_1 \leq u_1 = 0$. This means that 
$z \in P^{\overline{\{1\} \cup \sigma}}_{\mu} \cap \Z^m$, so that,
by Theorem \ref{thm:polytopeforstdpairs}, $z = 0$.
The rest of the assertions are trivial.
\end{proof}

\begin{theorem}
\label{thm:kernel}
$\ker(\partial_1) = \Span \big\{ x^u : u \in K \big\}$.
\end{theorem}

\begin{proof}

Let $\varphi$ be a (possibly logarithmic) solution of $H_A(\Av)$ 
such that $\partial_1 \varphi = 0$. The function $\varphi$ is a 
linear combination of canonical series with respect to $-e_1$.
Write $\varphi= \varphi_1+ \cdots + \varphi_l$  where each
$\varphi_i$ is a linear combination of canonical series
whose exponents differ by integer vectors, and the exponents
in $\varphi_i$ and $\varphi_j$ do not differ by integers if
$i \neq j$.

It is clear that $\partial_1 \varphi_i = 0$ for $ 1 \leq i \leq l$,
so we can reduce to the case when $\varphi$ is a 
linear combination of canonical series solutions whose exponents 
differ by integer vectors, and we assume this from now on.

Write $\varphi$ in the form of Observation \ref{keyobservation}
for some $\delta \in {\mathcal{S}}_{\max}$, so that
\[ \varphi = \varphi_{\delta}+ f_{\delta} \log(x)^{\delta} \]
where $\varphi_{\delta}$ contains only logarithmic terms that
less than or incomparable to $\delta$. Since $\partial_1 \varphi = 0$,
$f_{\delta}$ is a logarithm-free solution of $H_A(\Av)$ that 
is constant
with respect to $x_1$. Thus, $f_{\delta}$ is a linear combination of
logarithm-free canonical series whose corresponding
fake exponents differ by integer vectors, and have first
coordinate equal to zero. This means that $f_{\delta}$
is a linear combination of functions 
$\phi_{u^{(i)}}=x^{u^{(i)}}$ with $u^{(i)} \in K$, 
differing pairwise by integer vectors.

Now rewrite the function 
$\varphi$ in the form of Proposition \ref{propo:formoflogs},
that is $\varphi = \sum x^v g_v(\log(x))$. By the previous
reasoning, we can choose $u \in K$ such that $g_u \neq 0$.
We will show that $g_u$ is a constant. 

Let $(\partial^{\mu}, \{1 \} \cup \sigma)$ be a 
standard pair of $\ini_{-e_1}(I_A)$ such that 
$u_i=\mu_i$ for $i \not \in \{ 1 \} \cup \sigma$. By our choice of $\alpha$,
$u_i \in \Z$ implies $i \not \in \tau$, so that $\sigma \supseteq \tau =
\{ m+3, \dots , n\}$.

We first consider the case when $\sigma = \tau$. Pick $i \not \in \tau$.
By Theorem \ref{thm:polytopeforstdpairs} applied to 
$(\partial^{\mu}, \{1 \} \cup \tau)$, there exists $z \in \Z^m$
such that $(B\cdot z)_j \leq \mu_j$ for $j \not \in \{1,i\} \cup \tau$,
$(B\cdot z)_i > \mu_i$ and $(B\cdot z)_1 < 0$. 
Now choose a term $t^{\gamma}$ appearing with nonzero coefficient
in $g_u(t_1,\dots,t_n)$ such that $t_i$ appears to a maximal power
(among the terms in $g_u$).
Since 
$\partial_1 \varphi = 0$, 
\[ 0 = \partial^{(B\cdot z)_-} \varphi = \partial^{(B\cdot z)_+} \varphi. \]
Now, $\varphi$ contains a summand that is a nonzero multiple of
$x^u \log(x)^{\gamma}$. If $\gamma_i \neq 0$, then $\partial^{(B\cdot z)_+}
\varphi$ contains a term that is a nonzero multiple of
\begin{equation}
\label{eqn:nonzeroterm}
\frac{\big( \partial^{(B\cdot z)_+ -((B\cdot z)_i- u_i)e_i} x^u \big) 
\log(x)^{\gamma-e_i}}{x_i^{(B\cdot z)_i-u_i}} \;\; .
\end{equation}
By the construction of $z$, this term is nonzero. However, 
it cannot be cancelled with any other term from 
$\partial^{(B\cdot z)_+} \varphi = 0$, because such a term would have to 
come from $g_u$, and $\gamma_i$ was chosen maximal among those terms.
This contradiction implies that $\gamma_i = 0$. We conclude that 
$g_u$ is constant with respect to the $i$-th variable, and this
is true for all $i \not \in \{ 1 \} \cup \tau$. 

However, by Proposition \ref{propo:formoflogs} 
\[ g_u (t_1,\dots,t_n) = c_0 + c^{(1,1)}\cdot (t_1,\dots,t_n)+\cdots
+ \prod_{k=1}^l c^{(l,k)} \cdot (t_1,\dots,t_n) , \]
where $c_0 \in \C$, the vectors $c^{(j,k)}$ belong to the kernel of $A$,
and 
\[c^{(j,k)} \cdot (t_1,\dots ,t_n) = c^{(j,k)}_1t_1 + \cdots+
c^{(j,k)}_n t_n \; .\; \]
Since $g_u$ is constant with respect to the $i$-th variable for all $i \not
\in \{ 1\} \cup \tau$, $c^{(j,k)}_i = 0$ for $i \not \in \{1 \} \tau$,
$j=1,\dots l$, $1 \leq k \leq j$. By Lemma \ref{lemma:lind}, $m$
of the rows of $B$ indexed by $i \not \in \{ 1 \} \cup \tau$ are linearly
independent. This implies that all the vectors $c^{(j,k)}$ must
be zero, and we conclude that $g_u$ is constant.

Now we need to show that $g_u$ is constant in the case when
the standard pair corresponding to $u$ is $(\partial^{\mu},\{1\} \cup \sigma)$,
and $\sigma$ strictly contains $\tau$. In this case, 
$\sigma = \tau \cup \{ r \}$, for some $1 < r \leq m+2$.
Clearly, if $u_r \not \in \N$, the previous arguments still apply, because
no matter what $(B\cdot z)_r$ is, the term (\ref{eqn:nonzeroterm}) 
cannot vanish.
Thus, we may assume that $u_r \in \N$.

The polytope $P^{\overline{\{1\} \cup \sigma}}_{\mu}$ is a simplex, and
$0 \in P:= P^{\overline{\{1\} \cup \sigma}}_{\mu} \cap \{ z \in \R^m : 
(B\cdot z)_r \leq u_r \}$. This means that Lemma \ref{lemma:boundedregions}
applies, and so there is a set $T \subset \{ 1,\dots,\hat{r},\dots ,m+2 \}$
of cardinality $m$ such that, for each $i \in T$, the set obtained
from $P$ by reversing the inequality $(B\cdot z)_i \leq u_i$ 
is unbounded. Notice that, by Observation \ref{obs:changeofweight},
we may assume that $1 \not \in T$.

If $i \in T$, we can choose
$z \in \Z^m$ such that $(B\cdot z)_j \leq u_j$ for $j \not \in \{ 1, i \} \cup
\tau$, $(B\cdot z)_i > u_i$, and $(B\cdot z)_1 < 0$. Now we can
repeat what we did before, and conclude that $g_u$ is constant
with respect to the $i$-th variable, for all $i \in T$ .
Since the rows of $B$ indexed by $T$ are linearly independent, using
Proposition \ref{propo:formoflogs}, we conclude that $g_u$ is constant.

Now, remember that $g_u(\log(x))$ contains a term $\log(x)^{\delta}$
that is maximal among the logarithmic terms appearing in $\varphi$.
Since $g_u$ is constant, $\delta=0$, which implies that $\varphi$
is logarithm-free. But then it is clear that $\varphi$ belongs to 
$\Span \{ x^u : u \in K \}$.
\end{proof}

\section{The cokernel of the map $\partial_1$}
\label{sec:proofcokernel}

In this section we produce a subspace of
$\coker(\partial_1)$ whose dimension equals the dimension of
$\ker(\partial_1)$.

\begin{lemma}
\label{lemma:minimalsupport}
For each $u \in K$, the vector $u-e_1$ is a fake exponent with
minimal negative support of $H_A(\beta)$. 
\end{lemma}

\begin{proof}

Pick $u \in K$ corresponding to a standard pair 
$(\partial^{\mu}, \{ 1 \} \cup \sigma)$.
Clearly, $u-e_1$ is the fake exponent of $H_A(\beta)$
corresponding to this standard pair. Suppose that
$u-e_1$ does not have minimal negative support. Then we can find
$z \in \Z^m$ such that $u-e_1-B\cdot z$ has minimal negative
support strictly contained in $\nsupp(u-e_1)$. In particular, $z \neq 0$.
Notice that $1 \in \nsupp(u-e_1-B\cdot z)$, because otherwise
$z \in P^{\overline{\{1 \} \cup \sigma}}_{\mu} \cap \Z^m = \{ 0\}$.
This implies that the functional $-e_1$ is minimized
in the set $\{ B \cdot y : y \in \Z^m, \; \mbox{and} \; 
\nsupp(u-e_1-B\cdot z-B\cdot y) = \nsupp(u-e_1-B\cdot z) \}$. This 
minimum must be unique since $-e_1$ is a generic weight vector.
As a consequence, we may assume that $u-e_1-B\cdot z$ is a 
fake exponent of $H_A(\beta)$ with respect to the weight vector
$-e_1$. Then $u-B\cdot z$ is a fake exponent of
$H_A(\Av)$, it has minimal negative support, and $(u-B\cdot z)_1 < 0$.
The last two facts follow since $u$ has minimal negative support.
But now this contradicts Proposition \ref{propo:case1holds}.
\end{proof}

\begin{theorem}
\label{thm:cokernel}
The set $\Span \{\phi_{u-e_1} : u \in K \}$ does not intersect
the image of the map $\partial_1$.
\end{theorem}

\begin{proof}
Suppose there is a solution $\psi$ of $H_A(A\cdot v)$ such that
$\partial_1 \psi$ lies in $\Span \{\phi_{u-e_1} : u \in K \}$.
As in the proof of Theorem \ref{thm:kernel}, we may assume that
$\psi$ is a linear combination of canonical series whose
exponents differ by integer vectors.

We proceed as in the part of the proof of Theorem \ref{thm:kernel} where 
we show that the functions $\varphi_i$ are logarithm-free.
The first step is to write $\psi = \psi_{\delta}+f_{\delta} \log(x)^{\delta}$
for every $\delta \in {\mathcal{S}}_{\max}$ as in Observation
\ref{keyobservation}. The function $f_{\delta}$ belongs
to the kernel of $\partial_1$, so we may choose $u\in K$ such that 
$x^u$ appears with a nonzero coefficient in $f_{\delta}$.
Now we rewrite $\psi = \sum_v x^v g_v(\log(x))$ as in Proposition
\ref{propo:formoflogs}. We want to compute the polynomial $g_u$.

As in Theorem \ref{thm:kernel} we consider two cases, according
to whether the standard pair $(\partial^{\mu}, \{ 1 \} \cup \sigma)$
corresponding to the fake exponent $u$ is embedded or not. 

In the first case, $\sigma = \tau$, and $\nsupp(u)=\emptyset$. 
Choose $i \not \in \{1 \} \cup \tau$
and a vector $z \in \Z^m$ such that $(B\cdot z)_j \leq \mu_j$
for $j \not \in \{1,i\} \cup \tau$, $(B\cdot z)_i>\mu_i$, and
$(B\cdot z)_1<0$. We can do this by Theorem \ref{thm:polytopeforstdpairs}.

Now choose a term $t^{\gamma}$ appearing with nonzero coefficient
in the polynomial 
$g_u(t_1,\dots,t_n)$ such that $t_i$ appears to a maximal power
(among the terms in $g_u$), and consider the term (\ref{eqn:nonzeroterm})
from the proof of Theorem \ref{thm:kernel}:
\[ \frac{\big( \partial^{(B\cdot z)_+ -((B\cdot z)_i- u_i)e_i} x^u \big) 
\log(x)^{\gamma-e_i}}{x_i^{(B\cdot z)_i-u_i}} \;\;,\]
which appears with a nonzero coefficient in
$\partial^{(B\cdot z)_+} \psi=\partial^{(B\cdot z)_-} \psi$. 
Since $(B\cdot z)_1<0$, this function is a further derivative 
of $\partial_1 \psi \in \Span\{ \phi_{u-e_1}: u \in K\}$.

As $\nsupp(u)=\emptyset$, there are no fake exponents $v \in K$
such that $u-v \in \Z^m$.
Thus $\partial_1 \psi$ must be a multiple of $\phi_{u-e_1}$,
which has neither logarithmic terms, nor terms which
contain a strictly negative integer power of $x_i$. We conclude that
$\gamma_i=0$. The same argument as in Theorem \ref{thm:kernel}
now implies that $g_u$ is constant. This is a contradiction
because $x^{u-e_1}$ can only appear in $\partial_1 \psi$ if
$\psi$ has a term $x^u\log(x_1)$.

Now we need to consider the case when the standard pair
$(\partial^{\mu},\{ 1 \} \cup \sigma)$ is top-dimensional, that
is $\sigma = \tau \cup \{ r \}$, for some $2 \leq r \leq m+2$.
If $u_r \not \in \Z$, the same argument as above shows that
$g_u$ is constant. If $u_r \in \N$, we need to combine
the previous argument with the reasoning from Theorem
\ref{thm:kernel} to again conclude that $g_u$ is constant.
Finally if $u_r$ is a negative integer, we observe
that $\partial_1 \psi$ is a linear combination
of functions $\phi_{v-e_1}$, where the vectors 
$v$ belong to $K$ and differ with $u$ by an integer vector.
But then, for each $i \not \in \{ 1,r\} \cup \tau$,
a further derivative of $\partial_1 \psi$
has no logarithmic terms, or terms that contain strictly negative integer
powers of the variable $x_i$ and $x_r$. Unless $g_u$ is constant
with respect to the $i$-th variable, this
contradicts the fact that $\partial^{(B\cdot z)_+} \psi$
contains a nonzero multiple of the term (\ref{eqn:nonzeroterm}).
As before, we conclude that $g_u$ is constant, a contradiction.
\end{proof}

\section{Producing a rank jump}
\label{sec:proofjump}

We now analyze what solutions of $H_A(\Av)$ and $H_A(\beta)$
arise from the standard pair $(\partial^{\eta},\{ 1 \} \cup \tau)$.
Clearly, $\eta+\alpha$ is the fake exponent of $H_A(\Av)$
corresponding to this standard pair. By Lemma \ref{lemma:monomial},
the corresponding $A$-hypergeometric series is 
$\phi_{\eta+\alpha}=x^{\eta+\alpha}$.

For each $2 \leq i \leq m+2$, that is, for each $i \not \in \{ 1 \} \cup 
\tau$, we can choose $z^{(i)} \in \Z^m$ such that
$(B\cdot z^{(i)})_j \leq \eta_j$ for $j \not \in \{ 1,i\} \cup \tau$,
$(B\cdot z^{(i)})_i>\eta_i$ and $(B\cdot z^{(i)})_1 < 0$.
Let $z^{(1)}=0$, and
consider the vectors $\eta+\alpha-e_1-B\cdot z^{(i)}$. 
We have $\nsupp(\eta+\alpha-e_1-B\cdot z^{(i)})=\{ i \}$.
By Lemma \ref{lemma:minimalsupport}, 
$\eta+\alpha-e_1$ has minimal negative support.
It follows that all the $\eta+\alpha-e_1-B\cdot z^{(i)}$ have minimal 
negative support. Of course, this is not enough to guarantee that they 
all give rise
to logarithm-free solutions of $H_A(\beta)$. However, some of them do,
for instance $\eta+\alpha-e_1$, since it is a fake exponent.
The corresponding $A$-hypergeometric series is 
$\phi_1 := \phi_{\eta+\alpha-e_1}$.
The purpose of the next proposition is to show that there is 
at least one other solution of $H_A(\beta)$ arising this way.

\begin{proposition}
\label{propo:extrasols}
There exists at least one $i \not \in \{1 \} \cup \tau$ for which
we can find a vector $y^{(i)} \in \Z^m$ such that
$\eta+\alpha-e_1-B\cdot (z^{(i)}+y^{(i)})$ is a fake exponent with minimal
negative support equal to $\{ i \}$.
\end{proposition}

\begin{proof}

Let $N_i = \{ B\cdot z : z \in \Z^m \; \mbox{and} \;
\nsupp(u-e_1-B\cdot z)=\{ i \} \}$ for $i \not \in \{ 1 \} \cup \tau$.
If we can show that the linear functional $-e_1$ is minimized
uniquely in the set $U_i= \{ \eta+\alpha -B\cdot z : B\cdot z \in N_i\}$,
then this set will contain a fake exponent, which is what we want to
prove. Since the weight vector $-e_1$ is generic, it is enough to
see that $-e_1$ is minimized in $U_i$, or equivalently, that it
is maximized in $N_i$.

As in the proof of Observation \ref{obs:changeofweight} there are two cases.
Either the hyperplanes $\{ z : (B\cdot z)_i = 0\}$
are all distinct, for $i \not \in \tau$, or two of those coincide.

In the second case, we know that $\{ z : (B\cdot z)_1 = 0\}$
equals $\{ z : (B\cdot z)_r = 0\}$ for a certain $ 2 \leq r \leq m+2$.
But then, for each $i \not \in \{ 1,r \} \cup \tau$,
$B\cdot z \in N_i$ implies $(B\cdot z)_r \leq \eta_r$. Since the $r$-th
row of $B$ is a negative multiple of the first row,
the linear functional $-e_1$ is bounded above in $N_i$. Thus it
is maximized in $N_i$, and we obtain logarithm-free canonical solutions
$\phi_i$ of $H_A(\beta)$ corresponding to $\eta+\alpha-e_1 -B\cdot z^{(i)}$
for $i \not \in \{ 1, r \} \cup \tau$.

Now assume that the hyperplanes  $\{ z : (B\cdot z)_i = 0\}$
are pairwise different for $i \not \in \tau$, and remember
that $\langle \partial_i : i \not \in \{ 1 \} \cup \tau \rangle$ is an
embedded prime of $\ini_{-e_1}(I_A)$. Then there is a minimal prime
of this initial ideal of the form
$\langle \partial_i : i \not \in \{ 1,l \} \cup \tau \rangle$,
for some $l \not \in \{ 1 \} \cup \tau$. Since there must be standard pairs
corresponding to this minimal prime, we conclude that the set
obtained from $P^{\overline{\{ 1 \} \cup \tau}}_{\eta}$ by 
deleting the inequality $(B\cdot z)_l \leq \eta_l$ is a simplex.
Thus $N_l$ is bounded, so that $-e_1$ will be maximized
in $N_l$, and we obtain $\phi_l$, the corresponding solution
of $H_A(\beta)$. Notice that, by Lemma \ref{lemma:boundedregions},
the linear functional $-e_1$ attains no maximum in
$N_i$ for $i \neq 1,l$.
\end{proof}

The following theorem is the fundamental ingredient to produce a rank
jump.

\begin{theorem}
\label{thm:rankjump}
Let $s \geq 2$ be the dimension of the 
span of the functions constructed in Proposition \ref{propo:extrasols}
(this includes the function $\phi_1=\phi_{\eta+\alpha-e_1}$).
The intersection of this span with the image of $\partial_1$
has dimension at most $s-2$.
\end{theorem}

\begin{proof}

First notice that if $\psi$ is a logarithm-free solution
of $H_A(\Av)$, and $\partial_1 \psi$ lies in the span of
the functions from Proposition \ref{propo:extrasols},
then $\psi$ is a linear combination of canonical,
logarithm-free solutions of $H_A(\Av)$ that are either constant 
with respect to $x_1$, or whose fake exponents
differ by an integer vector with $\eta+\alpha$.
Since $\nsupp(\eta+\alpha)=\emptyset$, the only such fake exponent
is $\eta+\alpha$. Thus $\partial_1 \psi = 0$.

Now let $\psi$ be a logarithmic
solution of $H_A(\Av)$ such that $\partial_1 \psi$
lies in the span of the functions constructed in Proposition
\ref{propo:extrasols}. We may assume that
$\psi$ does not contain a summand $x^u$, $u \in K$.
As in the proof of Theorem \ref{thm:cokernel}, we can
write $\psi = \sum x^v g_v(\log(x))$, where $g_{\eta+\alpha}(\log(x))$
contains all the maximal logarithmic terms.

Pick $i \not \in \{1 \} \cup \tau$, and remember the vector
$z^{(i)} \in \Z^m$ from the paragraph before Proposition 
\ref{propo:extrasols}. An argument similar to that in the proof
of Theorem \ref{thm:cokernel} shows that, if $g_{\eta+\alpha}(\log(x))$
contains a power of $\log(x_i)$ greater than 1, then $\partial_1 \psi$
cannot be logarithm-free. This and Proposition \ref{propo:formoflogs}
imply that $g_{\eta+\alpha}(\log(x)) = c_1\log(x_1)+\cdots +c_n \log(x_n)$,
where the vector $(c_1,\dots ,c_n)$ lies in the kernel of $A$.
Thus,
\[ \psi = \psi_{\eta+\alpha}+
x^{\eta+\alpha}(c_1\log(x_1)+\cdots +c_n \log(x_n)) \; ,\]
where $\psi_{\eta+\alpha}$ is logarithm-free, and contains no 
terms $x^u$, $u \in K$. Notice that, if $\eta+\alpha-e_1+(B\cdot z)_i$
does not give rise to a solution of $H_A(\beta)$ as in Proposition
\ref{propo:extrasols}, then $c_i=0$. This follows from
the same arguments that proved Theorem \ref{thm:cokernel}.

We claim that, once the numbers $c_1,\dots ,c_n$ are fixed,
the function $\psi$ itself is fixed. This is because 
the difference between two such functions would be a logarithm-free
solution of $H_A(\Av)$ whose first derivative lies
in the span of the functions from Proposition \ref{propo:extrasols}.
By the first paragraph of this proof, this implies that
our difference lies in the kernel of $\partial_1$. But since it cannot
contain terms in $x^u$, $u \in K$, we conclude that this difference
must be zero.

Now we consider two cases, as in the proof of Proposition 
\ref{propo:extrasols}. In the case that a hyperplane
$\{z: (B\cdot z)_r=0\}$ coincides with $\{z: (B\cdot z)_1 = 0\}$,
we have $m+1$ linearly independent
functions $\phi_i$, $i \not \in \{r \} \cup \tau$
from Proposition \ref{propo:extrasols}. But if 
$\partial_1 \psi$ lies in the span of these functions,
we know that 
\[ \psi = \psi_{\eta+\alpha}+
x^{\eta+\alpha}(c_1\log(x_1)+\cdots +c_n \log(x_n)) \; ,\]
where $(c_1,\dots,c_n)$ belongs to the kernel of $A$,
and $c_r=0$. Since the dimension of the kernel of $A$ is $m$,
we conclude that the dimension of the space of such functions
$\psi$ is at most $m-1$. Thus, the intersection of
the image of $\partial_1$ and the span of the functions
$\phi_i$, $i \not \in \{r \} \cup \tau$ has dimension at most
$m-1=(m+1)-2=s-2$.

In the case when all the hyperplanes  $\{z: (B\cdot z)_i = 0\}$,
$i \not \in \tau$ are distinct, we have only two functions,
$\phi_1$ and $\phi_l$ from Proposition \ref{propo:extrasols}.
Thus, if $\partial_1 \psi$ lies in their span,
\[ \psi = \psi_{\eta+\alpha}+
x^{\eta+\alpha}(c_1\log(x_1)+\cdots +c_n \log(x_n)) \; ,\]
where $(c_1,\dots,c_n)$ belongs to the kernel of $A$,
and $c_i = 0$, for $i \not \in \{1,l\}\cup \tau$.
By Lemma \ref{lemma:lind}, this implies that
$(c_1,\dots,c_n)$ is the zero vector, so that $\psi$ vanishes.
Thus the intersection of $\Span\{\phi_1,\phi_l\}$ with the 
image of $\partial_1$ is $\{ 0 \}$, which has dimension
$s-2$.
\end{proof}

We are finally ready to prove Theorem \ref{thm:restricted}.

\begin{proof}[Proof of Theorem \ref{thm:restricted}]
In Theorem \ref{thm:cokernel}
we produced one function in the cokernel of $\partial_1$ for each 
function in the kernel of $\partial_1$. Furthermore, 
in Theorem \ref{thm:rankjump}, we produced at least
$2$ linearly independent functions in the cokernel of 
$\partial_1$ corresponding to $x^{\eta+\alpha}$. 
All these functions are clearly linearly independent
in the cokernel of $\partial_1$.
This means that:
\[ \dim(\coker(\partial_1)) \geq \dim(\ker(\partial_1))+1 \; . \]
Adding $\dim(\im(\partial_1))$ to both sides of this inequality
we obtain:
\[ \dim(\coker(\partial_1))+\dim(\im(\partial_1)) \geq
\dim(\ker(\partial_1))+\dim(\im(\partial_1)) +1 \; ,\]
or equivalently,
\[ \rank(H_A(\beta)) \geq \rank(H_A(A\cdot v)) + 1 \; .\]
Since $\rank(H_A(A\cdot v)) \geq \vol(A)$, it follows that
\[ \rank(H_A(\beta)) \geq \vol(A) + 1 \; . \]
\end{proof}


\section{An example}
\label{sec:example}

Even though generic toric ideals are common among all toric ideals,
it is hard in practice to construct examples of
generic configurations $A$ such that $\vol(A)$ is small. This is
a disadvantage when we want to perform rank computations using
computer algebra systems. In this section, we work out an example
where $I_A$ is {\em not} generic, although the techniques that proved
Theorem \ref{thm:restricted} are still successful.
Let 
\[ A= \left( \begin{array}{cccccc} 1 & 1 & 1 & 1 & 1 &1 \\
1 & 2 & 1 & 2 & 3 &0 \\
0 & 2 & 2 & 0 & 1 & 1
\end{array} \right) \; .\]
The toric ideal is:
\[ I_A = \left\langle \begin{array}{c}
\partial_3\partial_4-\partial_5\partial_6,
\partial_1\partial_2-\partial_5\partial_6,
\partial_3^3\partial_5-\partial_2^3\partial_6,
\partial_1\partial_3^2\partial_5-\partial_2^2\partial_4\partial_6 \\
\partial_1^2\partial_3\partial_5-\partial_2\partial_4^2\partial_6,
\partial_1^3\partial_5-\partial_4^3\partial_6,
\partial_2\partial_4^3-\partial_1^2\partial_5^2,
\partial_2^2\partial_4^2-\partial_1\partial_3\partial_5^2\\
\partial_2^3\partial_4-\partial_3^2\partial_5^2,
\partial_1\partial_3^3-\partial_2^2\partial_6^2,
\partial_1^2\partial_3^2-\partial_2\partial_4\partial_6^2,
\partial_1^3\partial_3-\partial_4^2\partial_6^2 \end{array}
\right\rangle \; . \]
Notice that $\vol(A)=8$, and that $I_A$ is not Cohen Macaulay.
We compute the initial ideal of $I_A$
with respect to $-e_1$: 
\[ \ini_{-e_1}(I_A) = \biggr\langle  \begin{array}{c}
\partial_5\partial_6,
\partial_3\partial_4-\partial_5\partial_6,
\partial_4^2\partial_6^2,
\partial_2\partial_4\partial_6^2,
\partial_2^2\partial_6^2,
\partial_4^3\partial_6,
\partial_2\partial_4^2\partial_6 \\
\partial_2^2\partial_4\partial_6,
\partial_3^3\partial_5-\partial_2^3\partial_6,
\partial_2\partial_4^3,
\partial_2^2\partial_4^2,
\partial_2^3\partial_4-\partial_3^2\partial_5^2 
\end{array} \biggr\rangle \; . \]
Since this ideal is not a monomial ideal, $I_A$ is not generic.
However, we can try to look at its embedded primes, in the spirit
of the generic case. Using the computer algebra system
{\tt Singular}, we can compute the set of associated primes of 
$\ini_{-e_1}(I_A)$. We obtain the following set:
\[ 
\left\{ \begin{array}{c}
\langle \partial_6,
\partial_5,
\partial_4 
\rangle,
\langle
\partial_6,
\partial_4,
\partial_3 \rangle,
\langle
\partial_6,
\partial_3,
\partial_2 \rangle,
\langle
\partial_5,
\partial_4,
\partial_2 \rangle \\
\langle
\partial_6,
\partial_4,
\partial_3,
\partial_2 \rangle,
\langle
\partial_5,
\partial_4,
\partial_3,
\partial_2 \rangle,
\langle
\partial_2,
\partial_3,
\partial_4,
\partial_5,
\partial_6 \rangle \end{array} \right\} \; . \]
Corresponding to the ideal $\langle \partial_5,\partial_4,\partial_3,
\partial_2 \rangle$, $\ini_{-e_1}(I_A)$ 
has a standard pair $(\partial_4,\{ 1, 6\})$, according to Altmann's
generalization of this notion (see \cite{chainprop}).
We can try to use the construction from Theorem \ref{thm:restricted}
applied to this standard pair. We obtain a line
\[ \{ A \cdot (-1,0,0,1,0,\alpha)=(\alpha,1,\alpha) : \alpha \in \C \} \]
of candidates for exceptional parameters.

We can check, using the computer algebra system {\tt Macaulay2},
that $\rank(H_A(0,1,0))=9$. The calculation, performed
in a {\tt Linux} computer with dual Pentium III-700/100 processors
and 512 MB of RAM, lasted about 30 minutes.

We would like to use the methods from the previous section to
prove that $(\alpha,1,\alpha)$ is an exceptional parameter.
The first difficulty we encounter is that we cannot use the
weight vector $-e_1$ to compute canonical series, since this 
weight vector is not generic.
To bypass this disadvantage, 
we consider the following initial ideal of $\ini_{-e_1}(I_A)$:
\[ \ini_w(I_A) = \biggr\langle \begin{array}{c}
\partial_5\partial_6,
\partial_4\partial_3,
\partial_4^2\partial_6^2,
\partial_2\partial_4\partial_6^2,
\partial_2^2\partial_6^2,
\partial_4^3\partial_6 \\
\partial_2\partial_4^2\partial_6,
\partial_2^2\partial_4\partial_6,
\partial_5\partial_3^3,
\partial_2\partial_4^3,
\partial_2^2\partial_4^2,
\partial_2^3\partial_4 \end{array} \biggr\rangle \; .\]
The set of standard pairs of this initial ideal is:
\[ \begin{array}{ll}
(\partial_3, \{1,2,3\}), & (1, \{1,4,5\}), \\
(1,\{1,2,3\}), & (\partial_4,\{1,6\}), \\
(\partial_2,\{1,3,6\}) , & (\partial_2\partial_4^2, \{1,5\}), \\
(1, \{ 1,3,6\}), & (\partial_2^2\partial_4,\{1,5\}), \\
(\partial_3^2, \{1,2,5\}), & (\partial_2\partial_4 , \{1,5\}), \\
(\partial_3, \{1,2,5\}), & (\partial_4^2\partial_6, \{ 1 \}), \\
(1,\{ 1,2,5\}), & (\partial_2\partial_4\partial_6, \{ 1 \}) .
\end{array} \]
Thus, the standard pair $(\partial_4,\{1,6\})$, which
we used to construct our candidate for exceptional parameter, is
also a standard pair on $\ini_w(I_A)$. We will use this (generic) weight
vector to compute canonical series solutions of
$H_A(A\cdot (0,0,0,1,0,\alpha))=H_A((1+\alpha,2,\alpha))$ 
and $H_A(A\cdot (-1,0,0,1,0,\alpha))=H_A((\alpha,1,\alpha))$.
As in the proof of Theorem \ref{thm:restricted}, we will use the
map $\partial_1$ between the solution spaces of 
$H_A((1+\alpha,2,\alpha))$ and $H_A((\alpha,1,\alpha))$.

If we chose $\alpha \neq 1,4,2/3$, so that the only embedded standard
pair that produces a fake exponent is $(\partial_4,\{1,6\})$,
the fake exponents of $H_A((1+\alpha,2,\alpha))$ with respect to $w$,
ordered by their corresponding standard pairs, are:
\[ \begin{array}{ll}
(\alpha/2+1/2,2-\alpha,3\alpha/2-5/2,0,0,1), 
& (3\alpha,0,0,1-3\alpha,\alpha,0), \\
(\alpha/2+1,1-\alpha,3\alpha/2-1,0,0,0), & 
(0,0,0,1,0,\alpha).\\
(3/2,0,1/2,0,0,\alpha-1), & \hspace{2cm} - \\
(1,1,-1,0,0,\alpha), &\hspace{2cm} - \\
(\alpha+2/3,\alpha-1/3,0,0,2/3-\alpha,0), &\hspace{2cm} - \\
(\alpha+1/3,\alpha-5/3,1,0,4/3-\alpha,0), &\hspace{2cm} - \\
(\alpha,\alpha-3,2,0,2-\alpha,0). &\hspace{2cm} - \end{array} \]
Now we choose $\alpha \neq -1,-2,-2/3,-1/3,0$, so that
none of the fake exponents of $H_A((1+\alpha,2,\alpha))$ 
have zero first coordinate except $(0,0,0,1,0,\alpha)$.
This implies that $\ker(\partial_1) =\Span \{\phi= x_4 x_6^{\alpha}\}$.

In order to prove that $\rank(H_A((\alpha,1,\alpha)))>\vol(A)=8$,
we need to produce at least two functions in the cokernel
of $\partial_1$. The first step is to make $\alpha$ generic
so that the only fake exponent of $H_A((1+\alpha,2,\alpha))$ that
differs with $(0,0,0,1,0,\alpha)$ by an integer vector is
$(1,1,-1,0,0,\alpha)$. Now, the vectors
$(-1,0,0,1,0,\alpha)$ and $(0,1,-1,0,0,\alpha)$ are fake exponents of 
$H_A((\alpha,1,\alpha))$ with minimal negative support.
Arguments similar to those in Section \ref{sec:proofjump} show that these
two functions are linearly independent elements of $\coker(\partial_1)$.
In conclusion,
\[\rank(H_A((\alpha,1,\alpha))) > 
\rank(H_A((1+\alpha,2,\alpha))) \geq \vol(A)\; . \]

Details on how to generalize this example can be found in \cite{mythesis}.


\section{The geometry of the exceptional set}
\label{sec:comprehensive-gb}

The first basic question about the exceptional set
is to determine exactly when it is nonempty. There are other
open problems in this area.
One of them is to determine 
whether or not this set is Zariski closed.
It is not hard to show, however, that it is Zariski constructible,
if we make use of comprehensive Gr\"obner bases.

\begin{proposition}
The exceptional set of a homogeneous matrix $A$ is Zariski
constructible, that is, it can be expressed as a finite Boolean combination
of Zariski closed sets.
\end{proposition}

\begin{proof}
We will prove our claim by presenting an algorithm to compute
${\mathcal{E}}(A)$ for a given homogeneous matrix $A$. This algorithm will 
rely on Gr\"obner basis computations in the Weyl algebra. 

In this proof we will think of $\beta$ as a parameter vector,
that is, a vector of indeterminates, instead of an element of $\C^d$. 
Thus, the $A$-hypergeometric system $H_A(\beta)$ will no longer be
an ideal in the Weyl algebra $D$, but an ideal in the {\em parametric}
Weyl algebra $D_{\C[\beta]}$, where the parameters $\beta_i$
commute with the variables $x_i$ and $\partial_i$.
Given a vector in $\C^d$, the {\em specialization map} corresponding
to this vector is the map from $D_{\C[\beta]}$ to $D$ that
replaces the parameters $\beta_i$ by the coordinates of our vector.
This allows us to introduce the concept of a comprehensive Gr\"obner basis
of a given left ideal $I$ in the parametric Weyl algebra.
Informally, a set $G$ is a {\em comprehensive left Gr\"obner basis} of $I$
if, for every vector in $\C^d$, the specialization of $G$ with respect to
that vector is a Gr\"obner basis of the corresponding specialization of $I$.

We are now ready to describe an algorithm to compute ${\mathcal{E}}(A)$.

\medskip

\noindent {\bf Input:} A homogeneous matrix $A$.

\noindent {\bf Output:} The exceptional set ${\mathcal{E}}(A)$.
\begin{enumerate} 
\item Compute a comprehensive left Gr\"obner
basis of $H_A(\beta)$
with respect to the weight vector 
whose first $n$ entries are zeros, and whose last $n$ entries
are ones. Call this basis $G$.
\item There are only finitely many
initial ideals $\ini_{(0,\dots,0,1,\dots,1)}(H_A(\beta))$
under specialization of the parameter $\beta$. These are ideals in the
polynomial ring $\C[x_1,\dots,x_n,s_1,\dots,s_n]$, call them
$I_1,\dots ,I_l$.
\item For each $1 \leq j \leq l$, the comprehensive
Gr\"obner basis $G$ yields finitely many polynomial conditions of the form
$g(\beta)=0$ or $g(\beta) \neq 0$ which characterize the (constructible)
subset $T_j$ of $\C^d$ such that if we specialize $\beta$ to an element
of $T_j$, the initial ideal of the specialized $A$-hypergeometric
system is $I_j$.
\item For $1 \leq j \leq l$ compute 
\[
r_j =  \dim_{\C(x_1,\dots ,x_n)}\left(\frac{\C(x_1,\dots ,x_n)[s_1,\dots,s_n]}
{\C(x_1,\dots ,x_n)[s_1,\dots,s_n] \cdot I_j} \right) .\]
\item ${\mathcal{E}}(A) = \cup_{j:r_j > \vol(A)} T_j.$
\end{enumerate}

By \cite[Formula 1.26]{SST}, the number $r_j$ produced
in the fourth step of the previous algorithm
equals the rank of any $D$-ideal whose initial ideal
with respect to $(0,\dots,0,1,\dots,1)$ is $I_j$. This 
justifies the description of ${\mathcal{E}}(A)$ from
the last step of the algorithm.
\end{proof}

Comprehensive Gr\"obner bases were introduced by 
Volker Weispfenning in \cite{compgb}.
This article deals only with the commutative case, and contains an explicit
algorithm for computing comprehensive Gr\"obner bases.
We refer to \cite{compsolv} for a proof of the existence of these objects
in the non commutative case.
Here, the authors argue that comprehensive G\"robner bases 
can be constructed, but they do not provide an explicit algorithm.
For this reason, it would be very desirable to have a method for computing
the exceptional set that did not require comprehensive Gr\"obner bases.

\noindent{\bf Acknowledgments:} I am very grateful to Bernd Sturmfels for 
many inspiring discussions, and unfailingly useful suggestions. 
I would also like to thank Diane Maclagan, for her comments on a previous
version of this article.


\def\cprime{$'$} \def\cprime{$'$}
\providecommand{\bysame}{\leavevmode\hbox to3em{\hrulefill}\thinspace}

\end{document}